\documentclass[11pt]{rMTA2010} % use larger type; default would be 10pt
\usepackage[utf8]{inputenc} % set input encoding (not needed with XeLaTeX)
\usepackage{graphicx} % support the \includegraphics command and options
\usepackage{booktabs} % for much better looking tables
\usepackage{array} % for better arrays (eg matrices) in maths
\usepackage{enumerate}
\usepackage{amsmath}
\usepackage{graphicx}
\usepackage{float}
\usepackage{subfig}
\usepackage{fancyhdr}
\paginas{1}{10} % Numeros de p�gina, la versi�n final la pone el editor; Page numbers, the Editor puts the final numbers

\begin{document}
\title{A New Method for the Analysis of Signals: The Square Wave Transform (SWT)} \author{
Osvaldo Skliar\thanks{Escuela de Inform\'atica, Universidad Nacional, Costa Rica. E-mail: oskliar@costarri-cense.cr}
\and Ricardo E. Monge\thanks{Escuela de Ciencias de la Computaci\'on e Inform\'atica, Universidad de Costa Rica, Costa Rica. E-mail: ricardo@mogap.net}
\and Guillermo Oviedo\thanks{Universidad Latina, San Pedro, Costa Rica. E-mail: oviedogmo@gmail.com}
\and Sherry Gapper\thanks{Universidad Nacional, Costa Rica. E-mail: sherry.gapper.morrow@una.cr}}

%\date{} % Activate to display a given date or no date (if empty),
         % otherwise the current date is printed

%\begin{document}
\maketitle
%$s^{-1}$
\begin{abstract}
The results obtained by analyzing signals with the Square Wave Method (SWM) introduced previously can be presented in the frequency domain clearly and precisely by using the Square Wave Transform (SWT) described here. As an example, the SWT is used to analyze a sequence of samples (that is, of measured values) taken from an electroencephalographic recording. A computational tool, available at www.appliedmathgroup.org/, has been developed and may be used to obtain automatically the SWTs of sequences of samples taken from registers of interest for biomedical purposes, such as those of an EEG or an ECG.
\end{abstract}

\KW  signal analysis, square wave method, square wave transform

\AMS 94A12, 65F99

\newpage
\begin{resumen}
%\hyphenpenalty=10000
Los resultados obtenidos al analizar se\~nales con el M\'etodo de las Ondas Cuadradas (Square Wave Method, SWM) ---previamente introducido--- pueden ser presentados en el dominio de la frecuencia de manera clara, precisa y concisa mediante el uso de la Transformada de las Ondas Cuadradas (Square Wave Transform, SWT). Se caracteriza la SWT y, como ejemplo, se la utiliza para analizar una secuencia de muestras (es decir, de valores medidos) tomadas de un registro electroencefalogr\'afico. En www.appliedmathgroup.org, se encuentra disponible un recurso computacional que posibilita obtener, de manera automatizada, las SWT de secuencias de muestras tomadas de registros de inter\'es biom\'edico como el EEG y el ECG, entre otros.
\end{resumen}
\PC an\'alisis de se\~nales, m\'etodo de las ondas cuadradas, transformada de las ondas cuadradas

\section{Introduction}

Consideration was previously given to the analysis of functions of one variable using the Square Wave Method (SWM) \cite{b1}. This method, which will be reviewed briefly in the following section, was generalized for functions of two variables and applied to the analysis of images \cite{b2}.

The objective of this article is to specify how the results obtained by analyzing signals with the SWM can be presented in the frequency domain clearly and concisely, using the mathematical process described below: the Square Wave Transform (SWT).

A preliminary version of this article was made available on arXiv \cite{b3a}.

\section{Brief review of the application of the SWM to the analysis of functions of one variable}

Given that this article is devoted to the analysis of signals, it will be considered that the independent variable is time ($t$).

Let $f(t)$ be a function of a variable $t$, which in a given interval $\Delta t$, satisfies the conditions of Dirichlet \cite{b3}: (1) In the interval $\Delta t$, the function $f(t)$ to be analyzed must have a finite number of relative maximums and minimums; (2) in that interval it also must have a finite number of points of discontinuity; and (3) for any instant of $\Delta t$, $f(t)$ must have a finite value. That function can then be approximated in that interval by means of a particular sum of trains of square waves. The use of the SWM makes it possible to specify these trains of square waves unambiguously.

Consider, for example, the function $f(t)$, as indicated below:

\begin{equation}
f(t)=3 \sin(2 \pi \cdot 5 \cdot t) + 4 \sin(2 \pi \cdot 7 \cdot t);\quad 0 \leq t \leq 1 \:\mathrm{s}
\end{equation}

\begin{figure}[H]
\centering
\includegraphics[width=3.5in]{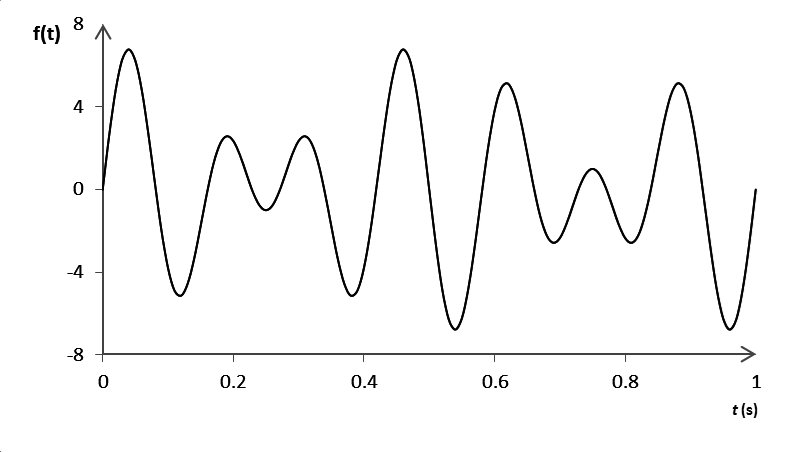}
\caption{$f(t)=3 \sin(2 \pi \cdot 5 \cdot t) + 4 \sin(2 \pi \cdot 7 \cdot t);\quad 0 \leq t \leq 1 \:\mathrm{s}$.}
\label{f1}
\end{figure}

In figure~\ref{f1}, $f(t)$ is shown for the interval specified in (1).

Note that the interval of $t$ ($\Delta t$), in which $f(t)$ will be analyzed, has a length of 1 second (1 s): $\Delta t = (1 - 0)\:\mathrm{s} = 1 \:\mathrm{s}$.

First, an explanation will be given about how to proceed if one wants to obtain an approximation to $f(t)$, in the interval $\Delta t$ specified in (1), composed of the sum of 10 trains of square waves. The interval $\Delta t$ is then divided into a number of sub-intervals -- of equal length -- which is the same as the number of trains of square waves. In this case, there will be 10 sub-intervals. The approximation to $f(t)$ to be obtained in interval $\Delta t$ will be the sum of 10 trains of square waves: $S_1$, $S_2$, $S_3$, \ldots, $S_9$, and $S_{10}$. The first of the trains of square waves will be referred to by $S_1$, the second by $S_2$, and so on.

Each of these trains of waves $S_i$, for $i = 1, 2, 3, \ldots, 9$ and $10$, will be characterized by a certain frequency $f_i$ (that is, the number of waves in the train of square waves considered which is contained in the unit of time), and a certain coefficient  $C_i$ whose absolute value is the amplitude of the corresponding train.

\begin{figure}[H]
\centering
\includegraphics[width=4in]{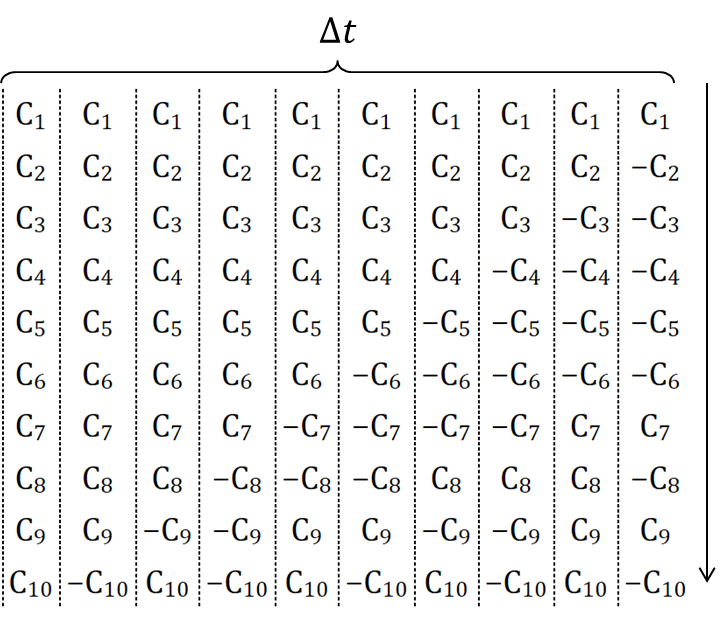}
\caption{How to apply the SWM to the analysis of the function represented in figure~\ref{f1}. (See indications in text.)}
\label{f2}
\end{figure}

For the case considered here, a description will be provided below of how the amplitudes corresponding to the different trains of square waves are determined (see figure~\ref{f2}).

The vertical arrow pointing down at the right of figure~\ref{f2} indicates how to add the terms corresponding to each of the 10 sub-intervals of $\Delta t$.  This procedure will make it possible to compute the values of the coefficients $C_1$, $C_2$, $C_3$, \ldots, $C_9$, and $C_{10}$ shown in figure 2. First, however, indications will be given about how to compute the frequencies $f_1$, $f_2$, $f_3$, \ldots, $f_9$, and $f_{10}$ corresponding respectively to the square wave trains $S_1$, $S_2$, $S_3$, \ldots, $S_9$, and $S_{10}$.

Each row in figure~\ref{f2} corresponds to the part of each of the trains of square waves $S_1, S_2, S_3, \dots, S_9$, and $S_{10}$ in interval $\Delta t$. In the first place, the structure of the last row in that figure corresponds to the part of the train of square waves $S_{10}$ in interval $\Delta t$. Each pair of consecutive coefficients \raisebox{-2pt}{\includegraphics[height=0.34cm]{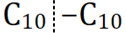}} corresponds to a square wave in $S_{10}$. Observe that in $\Delta t$ there are 5 of these pairs of elements (that is, there are 5 square waves in $\Delta t$). The frequency $f_{10}$, corresponding to $S_{10}$, is obtained by dividing, by $\Delta t$, the number of square waves occurring in $\Delta t$; thus, $f_{10} =\tfrac{5}{\Delta t}$. In the case discussed, $\Delta t = 1 \:\mathrm{s}$, so the value of $f_{10}$ is as follows: $f_{10} = 5\:\mathrm{s}^{-1}$. Of course, $\Delta t$ can be different from 1 s.  Suppose that we had taken $\Delta t=5\:\mathrm{s}$. The following value would have been obtained for $f_{10}$: $f_{10}=\tfrac{5}{5\:\mathrm{s}}=1\:\mathrm{s}^{-1}=$; and there would have been only one square wave in each time unit 1 s.

The next to the last row in figure~\ref{f2} corresponds to the part of the train of square waves $S_{9}$ in the interval $\Delta t$. Here the structure of each square wave is as follows: \raisebox{-2pt}{\includegraphics[height=0.34cm]{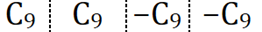}}. The length of the wave of each square wave corresponding to $S_{9}$ is double that of the wave corresponding to $S_{10}$. Note that each square wave corresponding to $S_{10}$ is included in 2 sub-intervals of $\Delta t$, whereas each square wave in $S_{9}$ is encompassed by 4 intervals of $\Delta t$. The value of $S_{9}$ is obtained by dividing, by $\Delta t$, the number of square waves $S_{9}$ in $\Delta t$: $f_9=\tfrac{2.5}{\Delta t}=2.5\:\mathrm{s}^{-1}$; in other words, in each unit of time 1 s, there are two and a half waves of $S_{9}$.) Observe that because the length of each square wave corresponding to $S_{9}$ is twice the length of each wave corresponding to $S_{10}$, the following result is to be expected: $f_9=\tfrac{1}{2}f_{10}$.
% \:\mathrm{s}
% \:\mathrm{s}^{-1}

The third row from the bottom in figure~\ref{f2} corresponds to the part of the train of square waves $S_{8}$ in the interval $\Delta t$. In this case, the structure of each square wave is: \raisebox{-2pt}{\includegraphics[height=0.34cm]{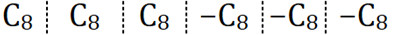}}. The length of each square wave in $S_{8}$ is three times that of each square wave in $S_{10}$. That is, each square wave in $S_{10}$ is encompassed by 2 sub-intervals of $\Delta t$, whereas each square wave of $C_{8}$ is encompassed by 6 sub-intervals of $\Delta t$. Of course, the value of $f_{8}$ is obtained by dividing, by $\Delta t$, the number of square waves corresponding to $S_{8}$ in the interval $\Delta t$: $f_8=\frac{\scriptscriptstyle \tfrac{10}{6}}{\tfrac{\Delta t}{1}}=\tfrac{1}{3}\cdot5\:\mathrm{s}^{-1}$.

Since the length of each square wave corresponding to $S_{8}$ is triple the length of each square wave corresponding to $S_{10}$, the validity of the following equality was foreseeable: $f_8=\tfrac{1}{3}f_{10}$.

In figure~\ref{f2}, it can be seen that the lengths of the square waves $S_7$, $S_6$, $S_5$, $S_4$, $S_3$, $S_2$, $S_1$ are respectively 4, 5, 6, 7, 8, 9, and 10 times longer than the square wave $S_{10}$.

Therefore, the following values are obtained for frequencies $f_1$, $f_2$, $f_3$, \ldots, $f_{10}$, corresponding respectively to $S_1, S_2, S_3, \dots, S_9$, and $S_{10}$:

\begin{alignat*}{3}
f_1 &= \frac{1}{10}\cdot f_{10} & &= \frac{1}{10}\cdot 5 \:\mathrm{s}^{-1} & &= 0.5000000\:\mathrm{s}^{-1} \\
f_2 &= \frac{1}{9}\cdot f_{10} & &= \frac{1}{9}\cdot 5 \:\mathrm{s}^{-1} & &=  0.5555556\:\mathrm{s}^{-1} \\
f_3 &= \frac{1}{8}\cdot f_{10} & &= \frac{1}{8}\cdot 5 \:\mathrm{s}^{-1} & &=  0.6250000\:\mathrm{s}^{-1} \\
f_4 &= \frac{1}{7}\cdot f_{10} & &= \frac{1}{7}\cdot 5 \:\mathrm{s}^{-1} & &=  0.7142857\:\mathrm{s}^{-1} \\
f_5 &= \frac{1}{6}\cdot f_{10} & &= \frac{1}{6}\cdot 5 \:\mathrm{s}^{-1} & &=  0.8333333\:\mathrm{s}^{-1} \\
f_6 &= \frac{1}{5}\cdot f_{10} & &= \frac{1}{5}\cdot 5 \:\mathrm{s}^{-1} & &=  1.0000000\:\mathrm{s}^{-1} \\
f_7 &= \frac{1}{4}\cdot f_{10} & &= \frac{1}{4}\cdot 5 \:\mathrm{s}^{-1} & &=  1.2500000\:\mathrm{s}^{-1} \\
f_8 &= \frac{1}{3}\cdot f_{10} & &= \frac{1}{3}\cdot 5 \:\mathrm{s}^{-1} & &=  1.6666667\:\mathrm{s}^{-1} \\
f_9 &= \frac{1}{2}\cdot f_{10} & &= \frac{1}{2}\cdot 5 \:\mathrm{s}^{-1} & &=  2.5000000\:\mathrm{s}^{-1} \\
f_{10} &= \frac{1}{1}\cdot f_{10} & &= \frac{1}{1}\cdot 5 \:\mathrm{s}^{-1} & &=  5.0000000\:\mathrm{s}^{-1} \\
\end{alignat*}

More concisely, these ten frequencies can be expressed as:
\begin{equation*}
f_i = \frac{1}{10-i+1}\cdot f_{10} = \frac{1}{10-i+1}\cdot 5\:\mathrm{s}^{-1};\quad i = 1, 2, 3, \dots, n
\end{equation*}

If the same approach is used for any $\Delta t$ expressed in seconds and any natural number $n$ of sub-intervals into which
the interval $\Delta t$ is divided, for the frequencies $f_1$, $f_2$, $f_3$, \ldots, $f_{n}$ corresponding respectively to the different trains of square waves $S_1, S_2, S_3, \dots, S_9$, and $S_{n}$, the following equation is obtained:
\begin{equation*}
f_i = \frac{1}{n-i+1}\cdot f_{n} = \frac{1}{n-i+1}\cdot \frac{\tfrac{n}{2}}{\Delta t} =  \frac{1}{n-i+1}\cdot \frac{n}{2\Delta t} ;\quad i = 1, 2, 3, \dots, n
\end{equation*}

How to compute the values of the coefficients $C_1$, $C_2$, $C_3$, \ldots, $C_9$, and $C_{10}$ shown in figure~\ref{f2} will be specified below.

First, the sum of all the coefficients in the first column of figure~\ref{f2} is made equal to the value of the function which one wants to approximate -- that is, (1) -- at the midpoint of the first of the ten sub-intervals into which the interval $\Delta t$ was divided. This value will be called $V_1$. Hence the following equation is obtained:
$C_1 + C_2 + C_3 + C_4 + C_5 + C_6 + C_7 + C_8 + C_9 + C_{10} = V_1$.

Second, the sum of all the coefficients in the second column of figure~\ref{f2} is made equal to the value of the function which one wants to approximate -- that is (1) -- at the midpoint of the second of the ten sub-intervals into which the interval  $\Delta t$ was divided. This value will be called $V_2$, and the following equation is obtained:
$C_1 + C_2 + C_3 + C_4 + C_5 + C_6 + C_7 + C_8 + C_9 - C_{10} = V_2$.

Third, the sum of all the coefficients in the third column of figure~\ref{f2} is made equal to the value of the function which one wants to approximate -- that is (1) -- at the midpoint of the third of the ten sub-intervals into which the interval $\Delta t$ was divided. This value will be called $V_3$, and the following equation is obtained:
$C_1 + C_2 + C_3 + C_4 + C_5 + C_6 + C_7 + C_8 - C_9 + C_{10} = V_3$.

The same is done for each of the remaining columns of coefficients in figure~\ref{f2}. Thus it is possible to obtain another seven equations which, together with the first three, constitute the following system of linear algebraic equations.

\begin{equation}
  \left.\begin{aligned}
C_1 + C_2 + C_3 + C_4 + C_5 + C_6 + C_7 + C_8 + C_9 + C_{10} &= V_1 \\
C_1 + C_2 + C_3 + C_4 + C_5 + C_6 + C_7 + C_8 + C_9 - C_{10} &= V_2 \\
C_1 + C_2 + C_3 + C_4 + C_5 + C_6 + C_7 + C_8 - C_9 + C_{10} &= V_3 \\
C_1 + C_2 + C_3 + C_4 + C_5 + C_6 + C_7 - C_8 - C_9 - C_{10} &= V_4 \\
C_1 + C_2 + C_3 + C_4 + C_5 + C_6 - C_7 - C_8 + C_9 + C_{10} &= V_5 \\
C_1 + C_2 + C_3 + C_4 + C_5 - C_6 - C_7 - C_8 + C_9 - C_{10} &= V_6 \\
C_1 + C_2 + C_3 + C_4 + C_5 - C_6 - C_7 + C_8 - C_9 + C_{10} &= V_7 \\
C_1 + C_2 + C_3 - C_4 - C_5 - C_6 - C_7 + C_8 - C_9 - C_{10} &= V_8 \\
C_1 + C_2 - C_3 - C_4 - C_5 - C_6 + C_7 + C_8 + C_9 + C_{10} &= V_9 \\
C_1 - C_2 - C_3 - C_4 - C_5 - C_6 + C_7 - C_8 + C_9 - C_{10} &= V_{10}
   \end{aligned}
  \qquad \right\}
\end{equation}

In the preceding system of linear algebraic equations (2), $V_1$, $V_2$, $V_3$, \ldots, $V_9$, and $V_{10}$ are the values for $f(t)$ as specified in (1) at the midpoints of the first, second, third, \ldots, ninth, and tenth sub-intervals, respectively, of the interval $\Delta t$ in which $f(t)$ is analyzed. It follows that the values $V_i$, for $i = 1,2, 3,\ldots, 9$ and $10$, can be computed given that $f(t)$ has been specified in (1). These values are:
\begin{align*}
& V_1 = 6.2360680 & &  V_6 = -6.2360680  \\
& V_2 = -1.7639320 & & V_7 = 1.7639320  \\
& V_3 = -1.0000000 & & V_8 = 1.0000000  \\
& V_4 = -1.7639320 & & V_9 = 1.7639320  \\
& V_5 = 6.2360680 & & V_{10} = -6.2360680
\end{align*}

The ten unknowns of the system of equations specified in (2) are $C_1$, $C_2$, $C_3$, \ldots, $C_9$, and $C_{10}$. $|C_i|$ refers to the amplitude of the train of square waves $S_i$, for $i = 1, 2, 3,\ldots,10$. The (constant) value of each positive square semi-wave of the train of square waves $S_i$ is $|C_i|$, and the (constant) value of each negative square semi-wave of that $S_i$ is $-|C_i|$.

The system of equations (2) was solved by using LAPACK \cite{b4}, and the following results were obtained for the unknowns:
\begin{align*}
& C_1 = -7.23607  & & C_6 = 2.23607\\
& C_2 = 3.61803 & & C_7 = 3.61803\\
& C_3 = 10.85410 & & C_8 = -3.61803\\
& C_4 = -3.61803 & & C_9 = 3.61803\\
& C_5 = -7.23607 & & C_{10} = 4.00000
\end{align*}

The trains of square waves $S_1, S_2, S_3, \dots, S_9$, and $S_{10}$ have been shown for interval $\Delta t$, in figures~\ref{f3a}, ~\ref{f3b}, ~\ref{f3c},\dots, ~\ref{f3i}, and ~\ref{f3j}, respectively.

\begin{figure}[H]
\centering
\subfloat[$S_1(t)$.]{
\includegraphics[width=4.5in]{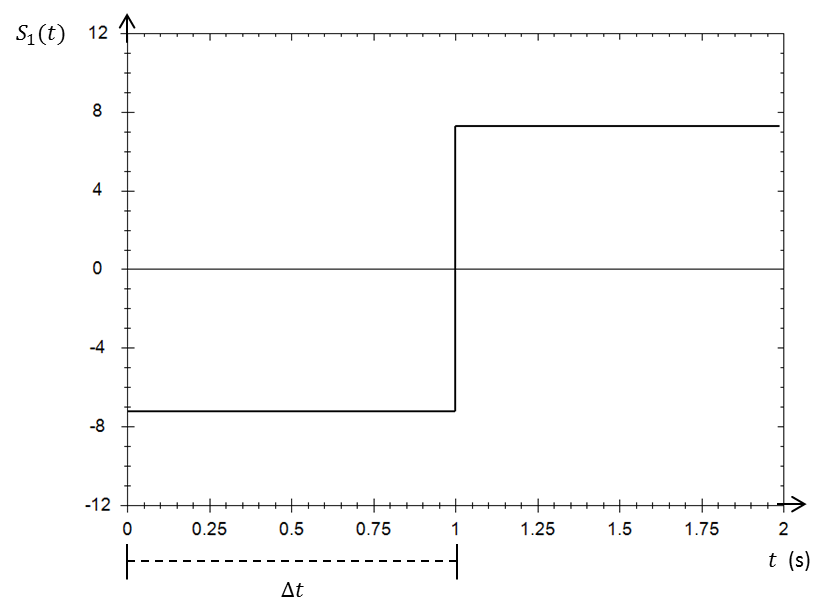}
\label{f3a}
}\qquad
\caption[]{}
\end{figure}
\begin{figure}
\ContinuedFloat
\subfloat[$S_2(t)$.]{
\includegraphics[width=4.5in]{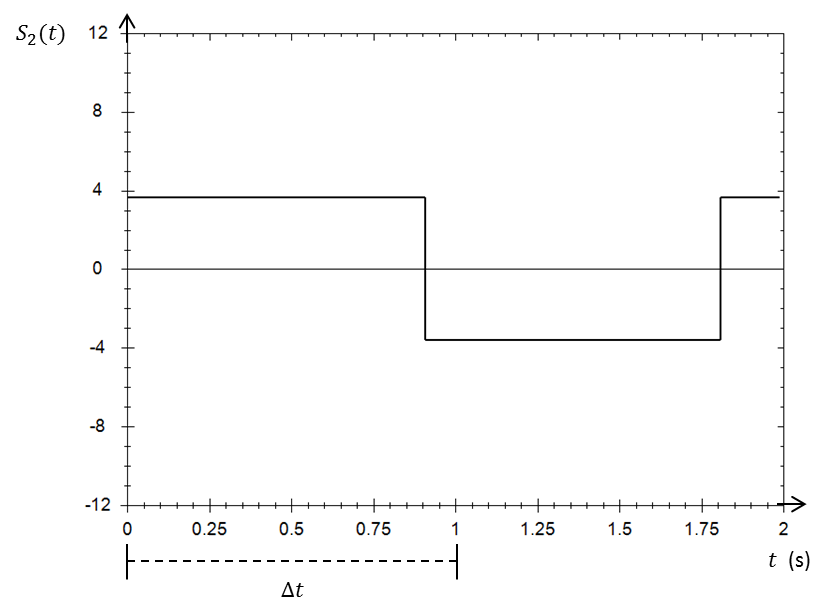}
\label{f3b}
}\qquad
\subfloat[$S_3(t)$.]{
\includegraphics[width=4.5in]{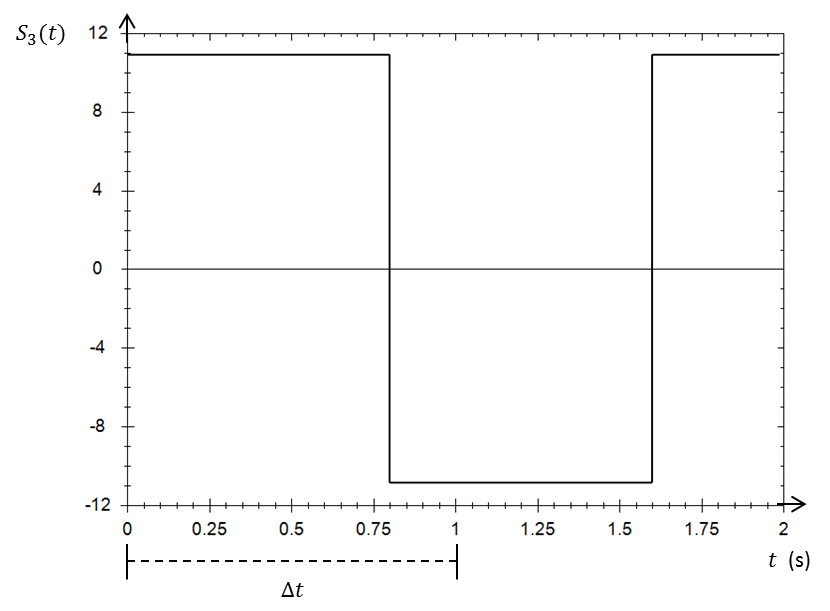}
\label{f3c}
}\qquad
\caption[]{}
\end{figure}
\begin{figure}
\ContinuedFloat
\subfloat[$S_4(t)$.]{
\includegraphics[width=4.5in]{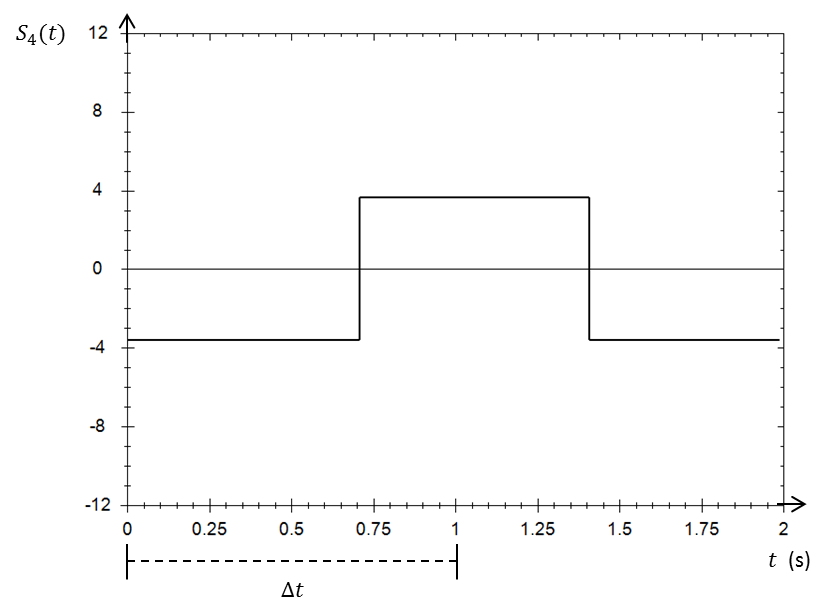}
\label{f3d}
}\qquad
\subfloat[$S_5(t)$.]{
\includegraphics[width=4.5in]{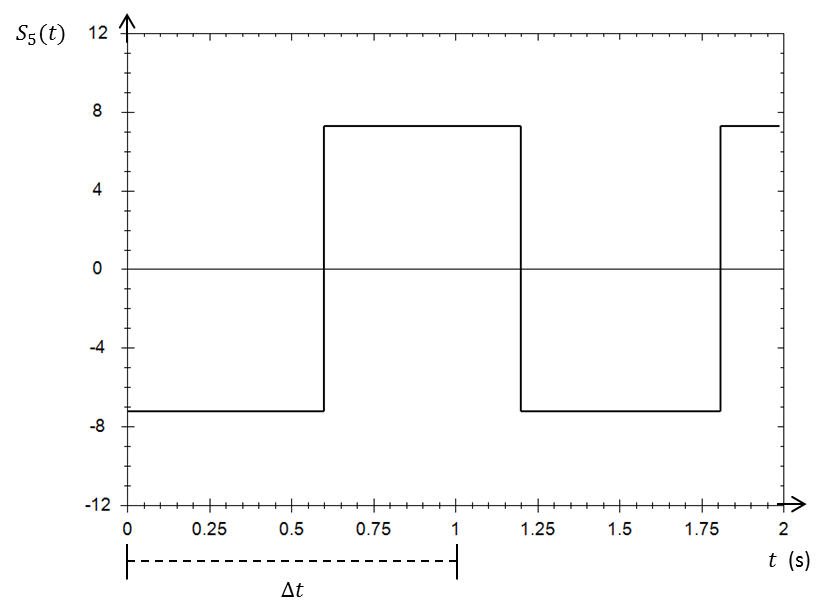}
\label{f3e}
}\qquad
\caption[]{}
\end{figure}
\begin{figure}
\ContinuedFloat
\subfloat[$S_6(t)$.]{
\includegraphics[width=4.5in]{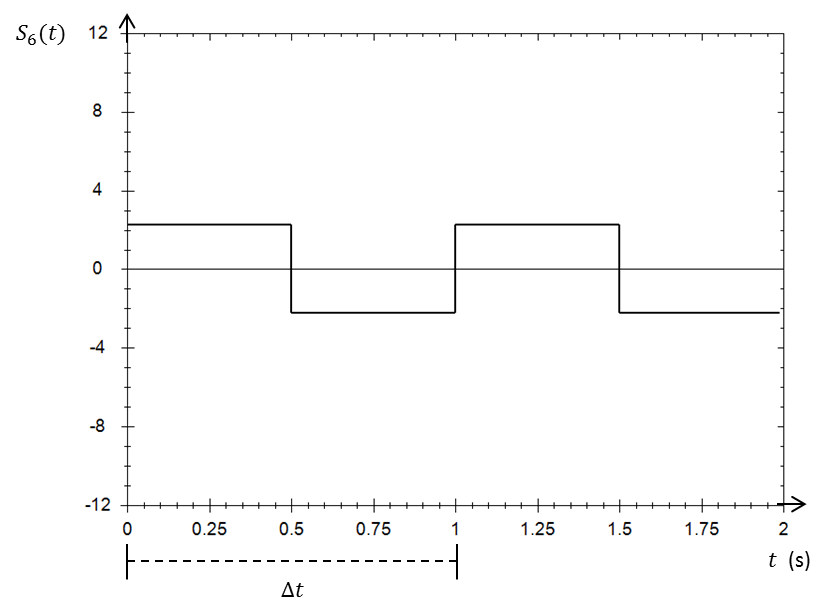}
\label{f3f}
}\qquad
\subfloat[$S_7(t)$.]{
\includegraphics[width=4.5in]{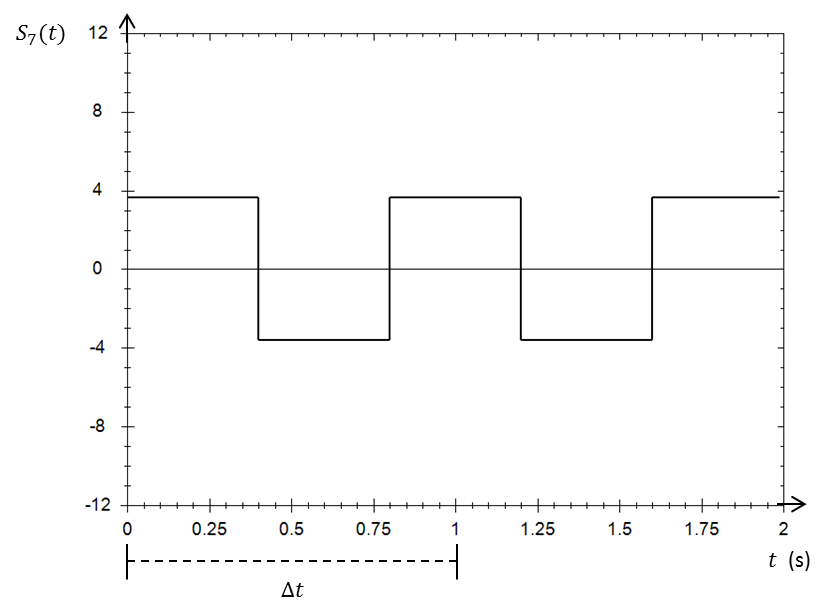}
\label{f3g}
}\qquad
\caption[]{}
\end{figure}
\begin{figure}
\ContinuedFloat
\subfloat[$S_8(t)$.]{
\includegraphics[width=4.5in]{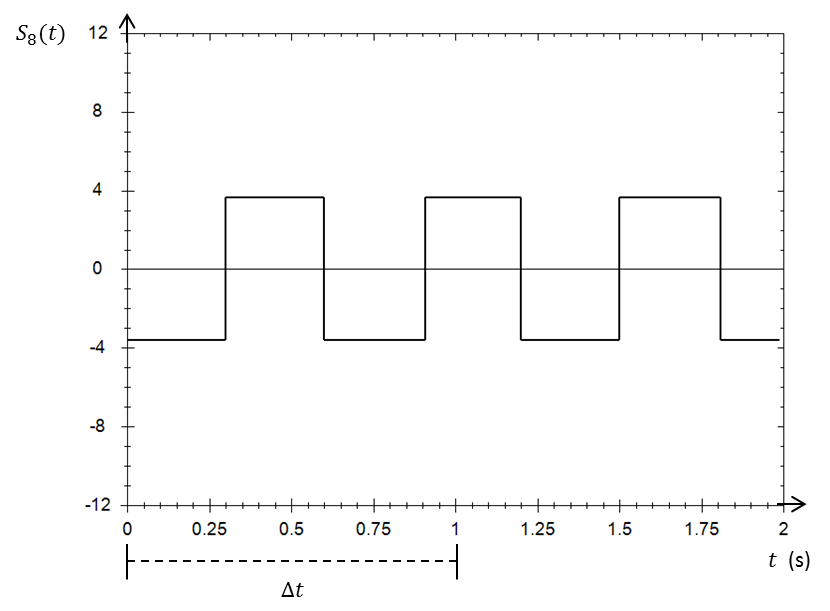}
\label{f3h}
}\qquad
\subfloat[$S_9(t)$.]{
\includegraphics[width=4.5in]{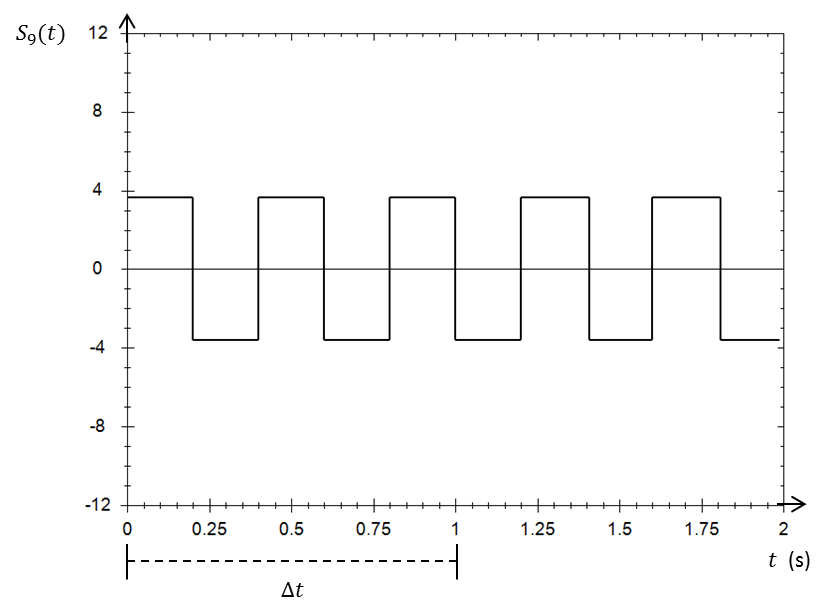}
\label{f3i}
}\qquad
\caption[]{}
\end{figure}
\begin{figure}
\ContinuedFloat
\subfloat[$S_{10}(t)$.]{
\includegraphics[width=5in]{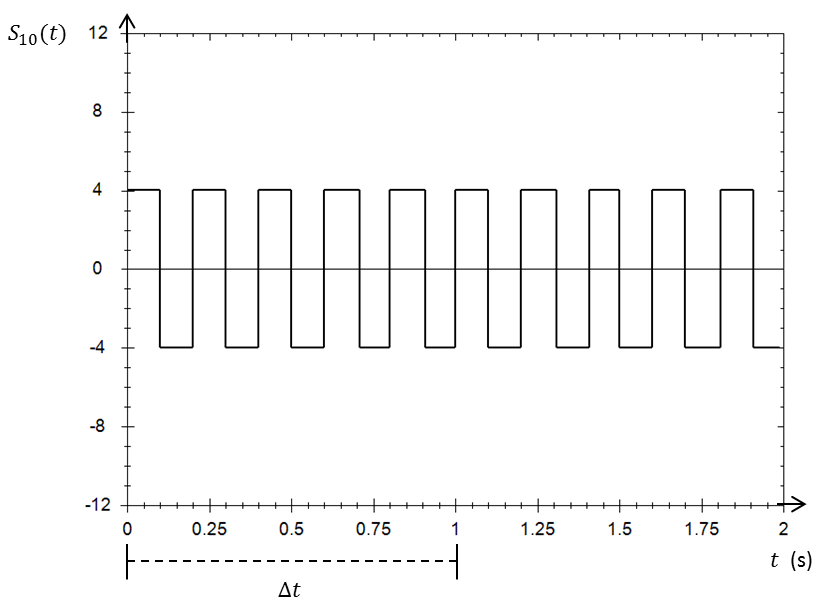}
\label{f3j}
}\qquad
\caption{Trains of square waves $S_1, S_2, S_3, \dots, S_9$, and $S_{10}$.}
\label{f3}
\end{figure}
\begin{figure}[H]
\centering
\includegraphics[width=4in]{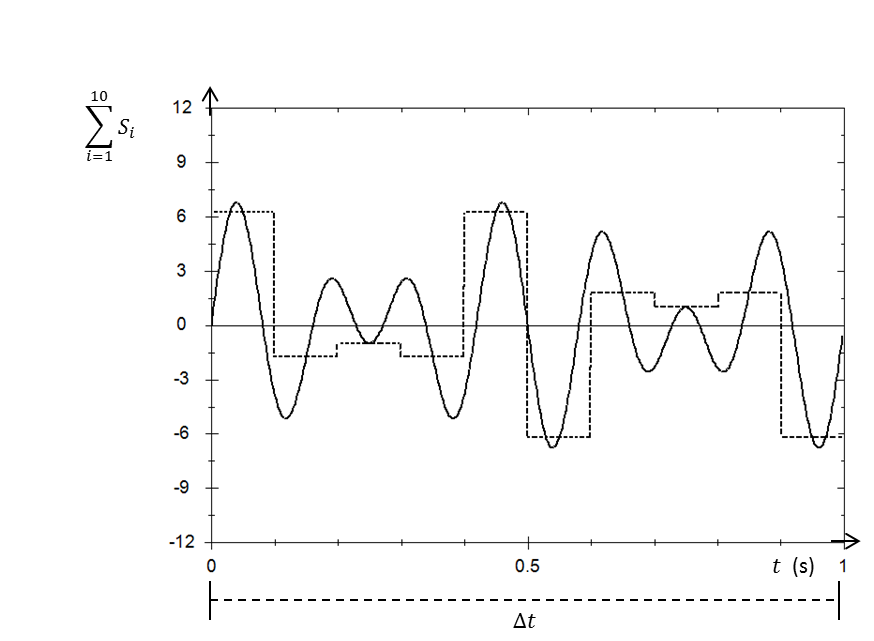}
\caption{The dashed line indicates the approximation to $f(t)$, specified in (1), by $\displaystyle\sum\limits_{i=1}^{10}S_ i$.}
\label{f4}
\end{figure}

\newpage
The approximation obtained for $f(t)$ (as specified in (1), in interval $\Delta t$, by adding the 10 trains of square waves) has been displayed in figure~\ref{f4}.

\begin{figure}[H]
\centering
\includegraphics[width=4.5in]{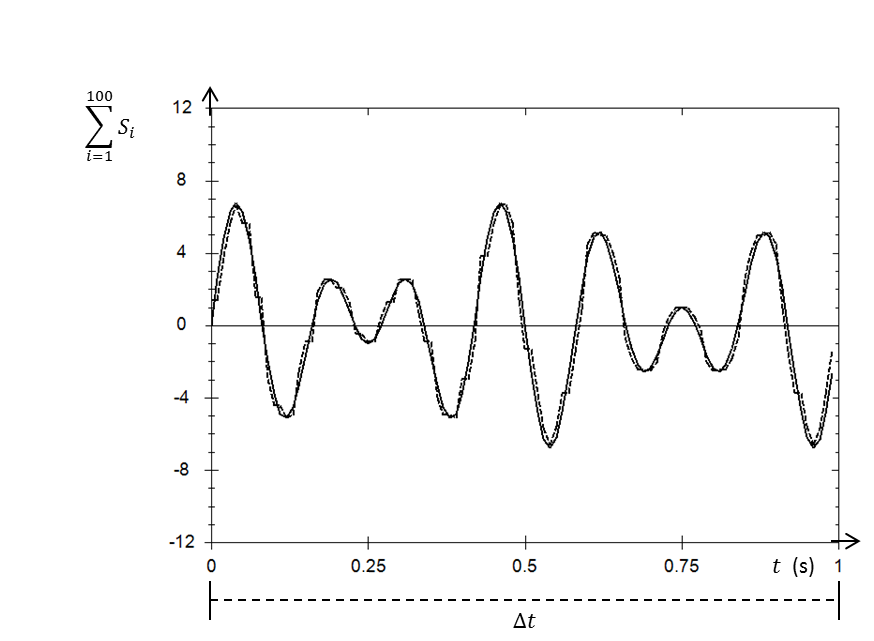}
\caption{The dashed line indicates the approximation to $f(t)$, specified in (1), by $\displaystyle\sum\limits_{i=1}^{100}S_ i$.}
\label{f5}
\end{figure}

If one wants to achieve a better approximation to $f(t)$, in interval $\Delta t$, by adding the trains of square waves, then $\Delta t$ should be divided into a larger number of equal sub-intervals. The larger the number of these sub-intervals, the better the approximation. Thus, for example, the approximation to $f(t)$ that can be achieved if $\Delta t$ is divided into 100 sub-intervals of equal length, is shown in figure~\ref{f5}. In this case, it is clear that 100 trains of square waves were added together.

The SWM cannot be considered a branch of Fourier analysis; that is, the trains of square waves $S_i$, for $i = 1, 2, 3,\dots, n$, do not make up a system of orthogonal functions.

\section{The Square Wave Transform (SWT) as a way of presenting the results of the analysis of $f(t)$ specified in (1)}

First, let us examine the results obtained when using the SWT to analyze the function $f(t)$ specified in (1), for the case specified above, in which the interval $\Delta t$ was divided into 10 equal sub-intervals. These results can be presented by a sequence of 10 dyads (ordered pairs) such that the first element of the first dyad is the frequency $f_1$ corresponding to $S_1$, and the second element of that first dyad is the coefficient $C_1$; the first element of the second dyad is the frequency $f_2$ corresponding to $S_2$, and the second element of that second dyad is the coefficient $C_2$; and so on successively, such that the first element of the tenth dyad is the frequency $f_{10}$ corresponding to $S_{10}$, and the second element of that tenth dyad is the coefficient $C_{10}$:
%7 then 5

\begin{equation*}
\begin{aligned}[c]
(f_{1}; C_{1}) &= (0.5000000; -7.23607)  \\
(f_{3}; C_{3}) &= (0.6250000; 10.85410)   \\
(f_{5}; C_{5}) &= (0.8333333; -7.23607)   \\
(f_{7}; C_{7}) &= (1.2500000; 3.61803)\\
(f_{9}; C_{9}) &= (2.5000000; 3.61803)\\
\end{aligned}
\qquad
\begin{aligned}[c]
(f_{2}; C_{2})&=    (0.5555556; 3.61803) \\
(f_{4}; C_{4}) &=   (0.7142857; -3.61803)   \\
(f_{6}; C_{6}) &=   (1.0000000; 2.23607)  \\
(f_{8}; C_{8}) &=   (1.6666667; -3.61803)\\
(f_{10}; C_{10}) &= (5.0000000; 4.00000) \\
\end{aligned}
\end{equation*}

\begin{figure}[H]
\centering
\includegraphics[width=4in]{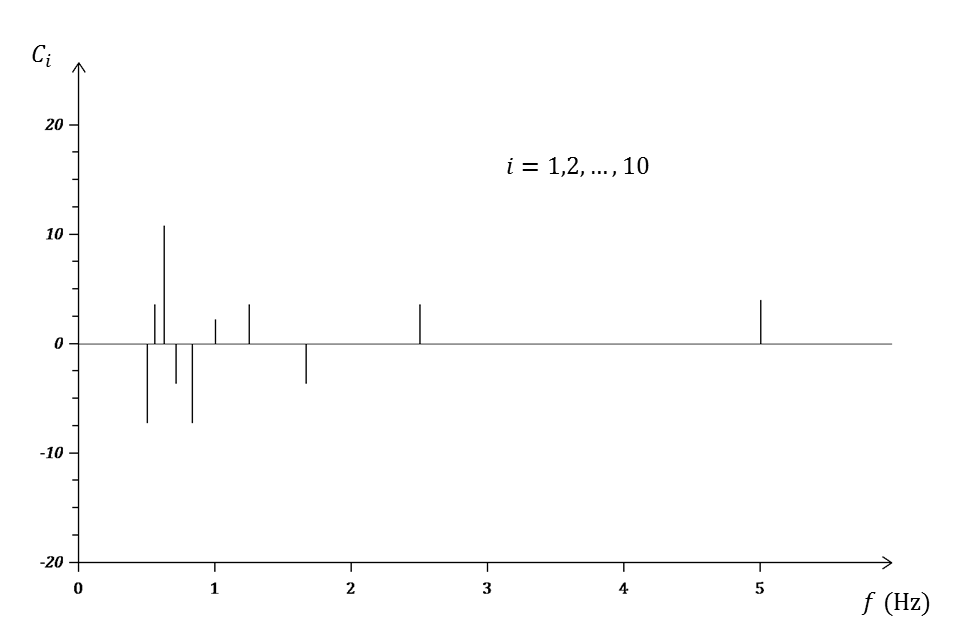}
\caption{SWT of the approximation to $f(t)$, specified in (1), obtained by dividing the interval $\Delta t$ into 10 sub-intervals.}
\label{f6}
\end{figure}

This approximation of the function $f(t)$ specified in (1) can be expressed in the frequency domain. To achieve this objective, for each of the frequencies considered, $f_1$, $f_2$, $f_3$, \ldots, $f_{10}$, the corresponding coefficients $C_1$, $C_2$, $C_3$, \ldots, $C_{10}$ must be indicated. The expression in the frequency domain of this approximation to $f(t)$ will be called the Square Wave Transform (SWT) of the approximation to that $f(t)$. This SWT is displayed in figure~\ref{f6}.

Of course, the SWTs corresponding to numbers as large as desired of equal sub-intervals into which $\Delta t$ is divided can be obtained for the $f(t)$ specified in (1), or for any other function of the time which, in a particular interval $\Delta t$, satisfies the conditions of Dirichlet.

Here, the symbol $N_{s}$ will be used to refer to the number of equal sub-intervals into which $\Delta t$ is divided. In figure~\ref{f7}, the SWTs of the approximations (to the $f(t)$ specified in (1)) obtained are shown for the following values of $N_{s}$: 100, 200, and 400.

\begin{figure}
\centering
\subfloat[SWT corresponding to $N_s=100$.]{
\includegraphics[width=4.5in]{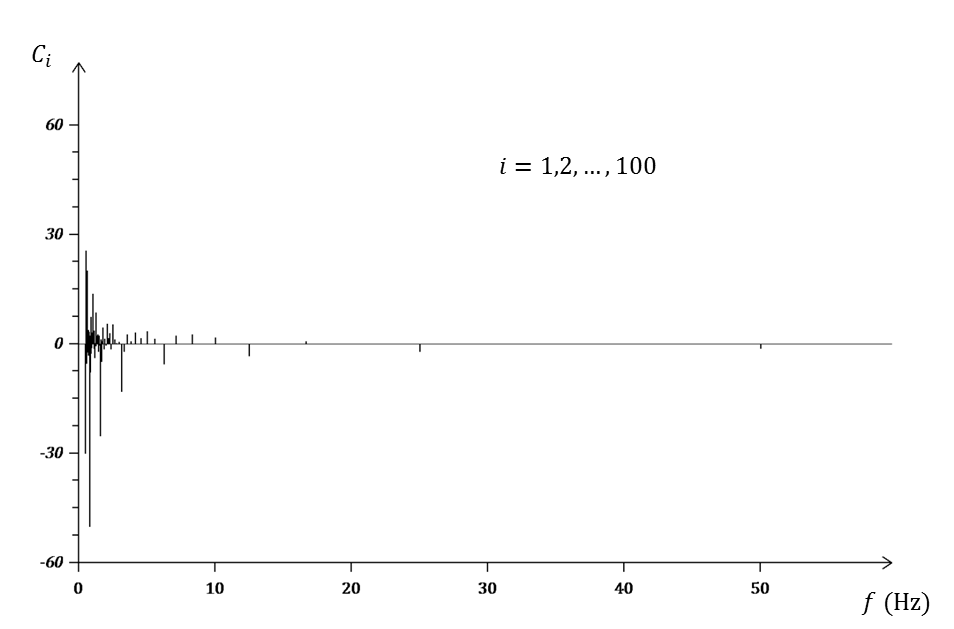}
\label{f7a}
}\qquad
\subfloat[SWT corresponding to $N_s=200$.]{
\includegraphics[width=4.5in]{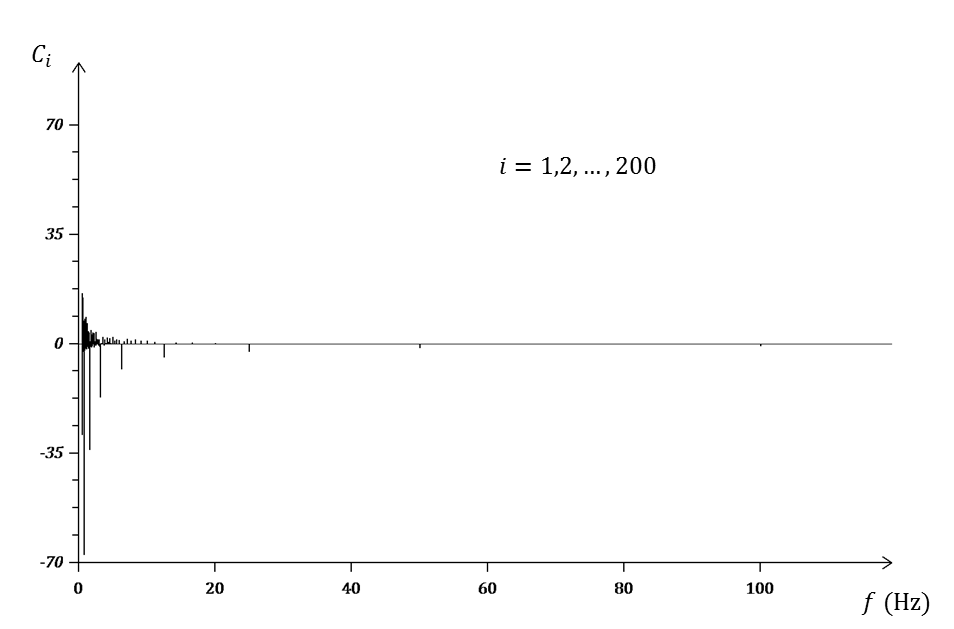}
\label{f7b}
}\qquad
\caption[]{}
\end{figure}

\begin{figure}
\ContinuedFloat
\subfloat[SWT corresponding to $N_s=400$.]{
\includegraphics[width=4.5in]{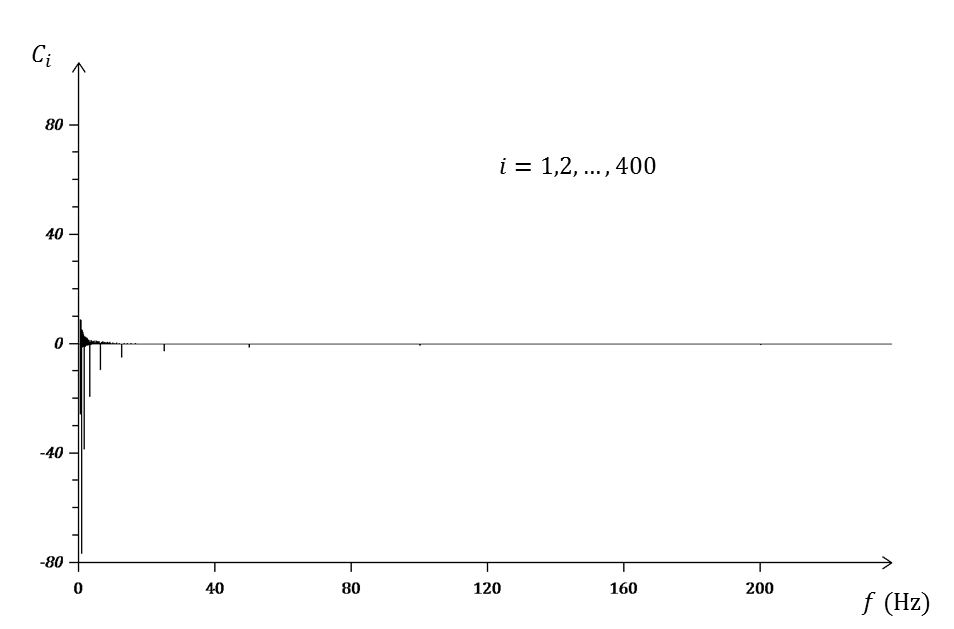}
\label{f7c}
}\qquad
\caption{The SWTs obtained of the approximations to the function $f(t)$ (specified in (1)), have been displayed in 7a, 7b, and 7c, for $N_{s}=100$, $N_{s}=200$, and $N_{s}=400$, respectively.}
\label{f7}
\end{figure}
\newpage
Note that in the three cases displayed in figure~\ref{f7}, different scales were used for the axes of the abscissas. The same scale will be used for the axes in figure~\ref{f10}.

\section{The SWT as a tool for the analysis of an electroencephalographic signal}

The SWT can be used for the analysis of signals of biomedical interest, such as those of electrocardiograms (ECG), electroencephalograms (EEG), electromiograms (EMG), etc.

Suppose that one has a sequence of 10 values of an electrophysiological signal, such as an electroencephalographic recording. To obtain the SWT corresponding to that sequence, the sequence of values is treated the same as was treated, with the same objective (that of obtaining the SWT) the sequence of values $V_1$, $V_2$, $V_3$, \ldots, $V_{10}$, in the system of algebraic equations (2). Generally, if one wants to obtain the SWT corresponding to a sequence $V_1$, $V_2$, $V_3$, \ldots, $V_{N}$ of measured values from that recording, that sequence of values is treated the same as the sequence of values  $V_1$, $V_2$, $V_3$, \ldots, $V_{N_s}$ is treated, with $N_{s}=N$. Let us recall that $V_1$, $V_2$, $V_3$, \ldots, $V_{N_s}$ is a sequence of values of a function which are computed at the midpoints of the $N_{s}$ equal sub-intervals of interval $\Delta t$ for which the function is characterized analytically.

\begin{figure}[H]
\centering
\includegraphics[width=4.5in]{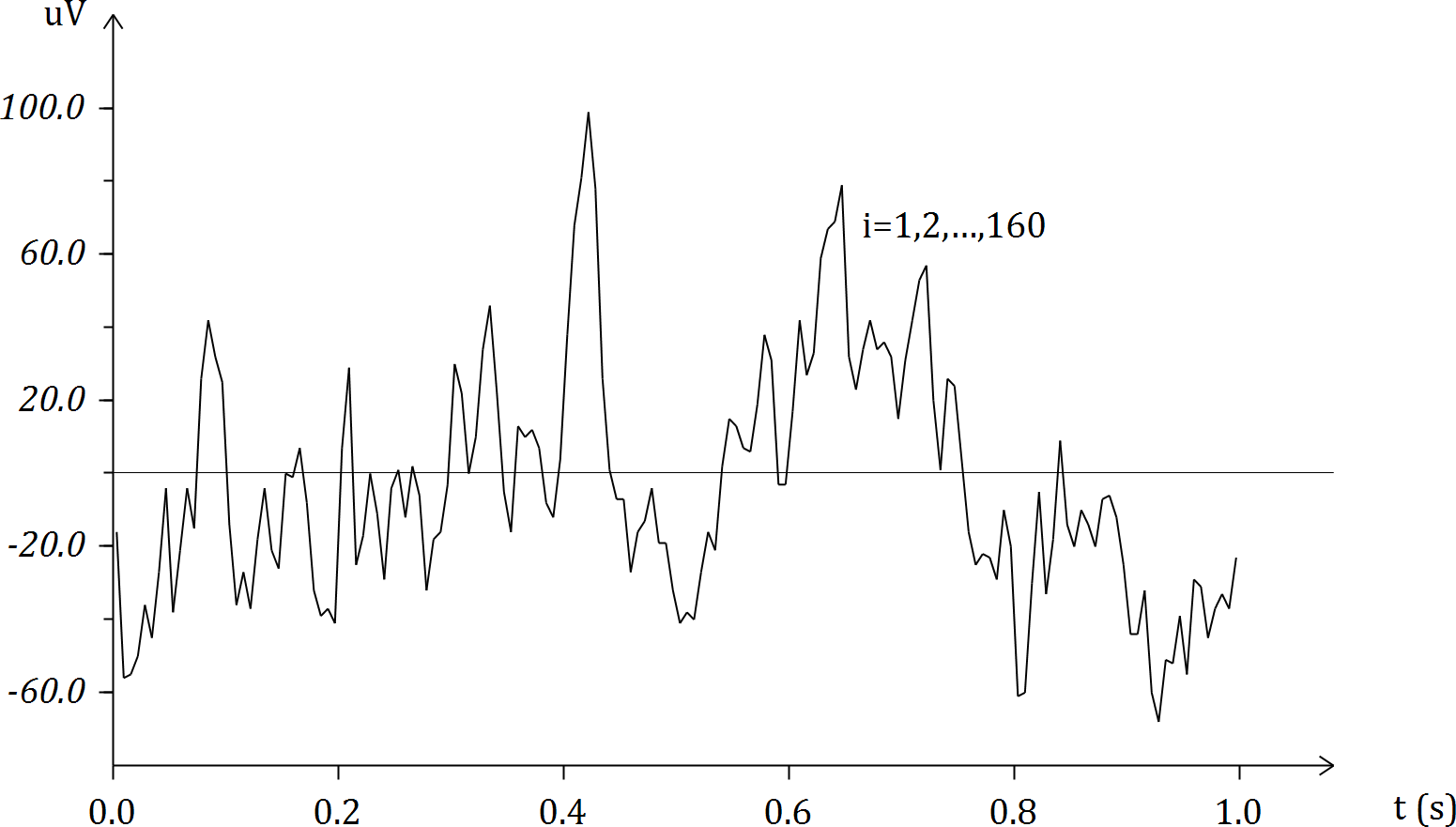}
\caption{Excerpt from an EEG, 160 Hz.}
\label{f8}
\end{figure}

A sequence of 160 ``samples'' (i.e., measured values) from an electroencephalographic recording is displayed in figure~\ref{f8}. The data were taken from the EEG Motor Movement/ Imagery Dataset (tagged MMIDB) in PhysioBank \cite{b5}. That recording corresponds to FC5 of run 01 of Subject S001, with a ``sampling'' frequency of 160 Hz.

The sequence of voltage measurements specified in microvolts ($\mathrm{\mu}V$), shown in figure~\ref{f8},  is as follows:

\vspace{4mm}

\begin{tabular}{ l l l l l }
$V_{1}=-16$ & $V_{33}=7$ & $V_{65}=38$ & $V_{97}=17$ & $V_{129}=-61$ \\
$V_{2}=-56$ & $V_{34}=29$ & $V_{66}=68$ & $V_{98}=42$ & $V_{130}=-60$ \\
$V_{3}=-55$ & $V_{35}=-25$ & $V_{67}=81$ & $V_{99}=27$ & $V_{131}=-30$ \\
$V_{4}=-50$ & $V_{36}=-17$ & $V_{68}=99$ & $V_{100}=33$ & $V_{132}=-5$ \\
$V_{5}=-36$ & $V_{37}=0$ & $V_{69}=78$ & $V_{101}=59$ & $V_{133}=-33$ \\
$V_{6}=-45$ & $V_{38}=-11$ & $V_{70}=26$ & $V_{102}=67$ & $V_{134}=-18$ \\
$V_{7}=-27$ & $V_{39}=-29$ & $V_{71}=1$ & $V_{103}=69$ & $V_{135}=9$ \\
$V_{8}=-4$ & $V_{40}=-4$ & $V_{72}=-7$ & $V_{104}=79$ & $V_{136}=-14$ \\
$V_{9}=-38$ & $V_{41}=1$ & $V_{73}=-7$ & $V_{105}=32$ & $V_{137}=-20$ \\
$V_{10}=-21$ & $V_{42}=-12$ & $V_{74}=-27$ & $V_{106}=23$ & $V_{138}=-10$ \\
$V_{11}=-4$ & $V_{43}=2$ & $V_{75}=-16$ & $V_{107}=34$ & $V_{139}=-14$ \\
$V_{12}=-15$ & $V_{44}=-6$ & $V_{76}=-13$ & $V_{108}=42$ & $V_{140}=-20$ \\
$V_{13}=26$ & $V_{45}=-32$ & $V_{77}=-4$ & $V_{109}=34$ & $V_{141}=-7$ \\
$V_{14}=42$ & $V_{46}=-18$ & $V_{78}=-19$ & $V_{110}=36$ & $V_{142}=-6$ \\
$V_{15}=32$ & $V_{47}=-16$ & $V_{79}=-19$ & $V_{111}=32$ & $V_{143}=-12$ \\
$V_{16}=25$ & $V_{48}=-3$ & $V_{80}=-32$ & $V_{112}=15$ & $V_{144}=-25$ \\
$V_{17}=-14$ & $V_{49}=30$ & $V_{81}=-41$ & $V_{113}=31$ & $V_{145}=-44$ \\
$V_{18}=-36$ & $V_{50}=22$ & $V_{82}=-38$ & $V_{114}=42$ & $V_{146}=-44$ \\
$V_{19}=-27$ & $V_{51}=0$ & $V_{83}=-40$ & $V_{115}=53$ & $V_{147}=-32$ \\
$V_{20}=-37$ & $V_{52}=10$ & $V_{84}=-27$ & $V_{116}=57$ & $V_{148}=-60$ \\
$V_{21}=-18$ & $V_{53}=34$ & $V_{85}=-16$ & $V_{117}=20$ & $V_{149}=-68$ \\
$V_{22}=-4$ & $V_{54}=46$ & $V_{86}=-21$ & $V_{118}=1$ & $V_{150}=-51$ \\
$V_{23}=-21$ & $V_{55}=22$ & $V_{87}=2$ & $V_{119}=26$ & $V_{151}=-52$ \\
$V_{24}=-26$ & $V_{56}=-5$ & $V_{88}=15$ & $V_{120}=24$ & $V_{152}=-39$ \\
$V_{25}=0$ & $V_{57}=-16$ & $V_{89}=13$ & $V_{121}=4$ & $V_{153}=-55$ \\
$V_{26}=-1$ & $V_{58}=13$ & $V_{90}=7$ & $V_{122}=-16$ & $V_{154}=-29$ \\
$V_{27}=7$ & $V_{59}=10$ & $V_{91}=6$ & $V_{123}=-25$ & $V_{155}=-31$ \\
$V_{28}=-8$ & $V_{60}=12$ & $V_{92}=19$ & $V_{124}=-22$ & $V_{156}=-45$ \\
$V_{29}=-32$ & $V_{61}=7$ & $V_{93}=38$ & $V_{125}=-23$ & $V_{157}=-37$ \\
$V_{30}=-39$ & $V_{62}=-8$ & $V_{94}=31$ & $V_{126}=-29$ & $V_{158}=-33$ \\
$V_{31}=-37$ & $V_{63}=-12$ & $V_{95}=-3$ & $V_{127}=-10$ & $V_{159}=-37$ \\
$V_{32}=-41$ & $V_{64}=4$ & $V_{96}=-3$ & $V_{128}=-20$ & $V_{160}=-23$ \\
\end{tabular}

\vspace{8mm}

\begin{figure}[H]
\centering
\includegraphics[width=4.5in]{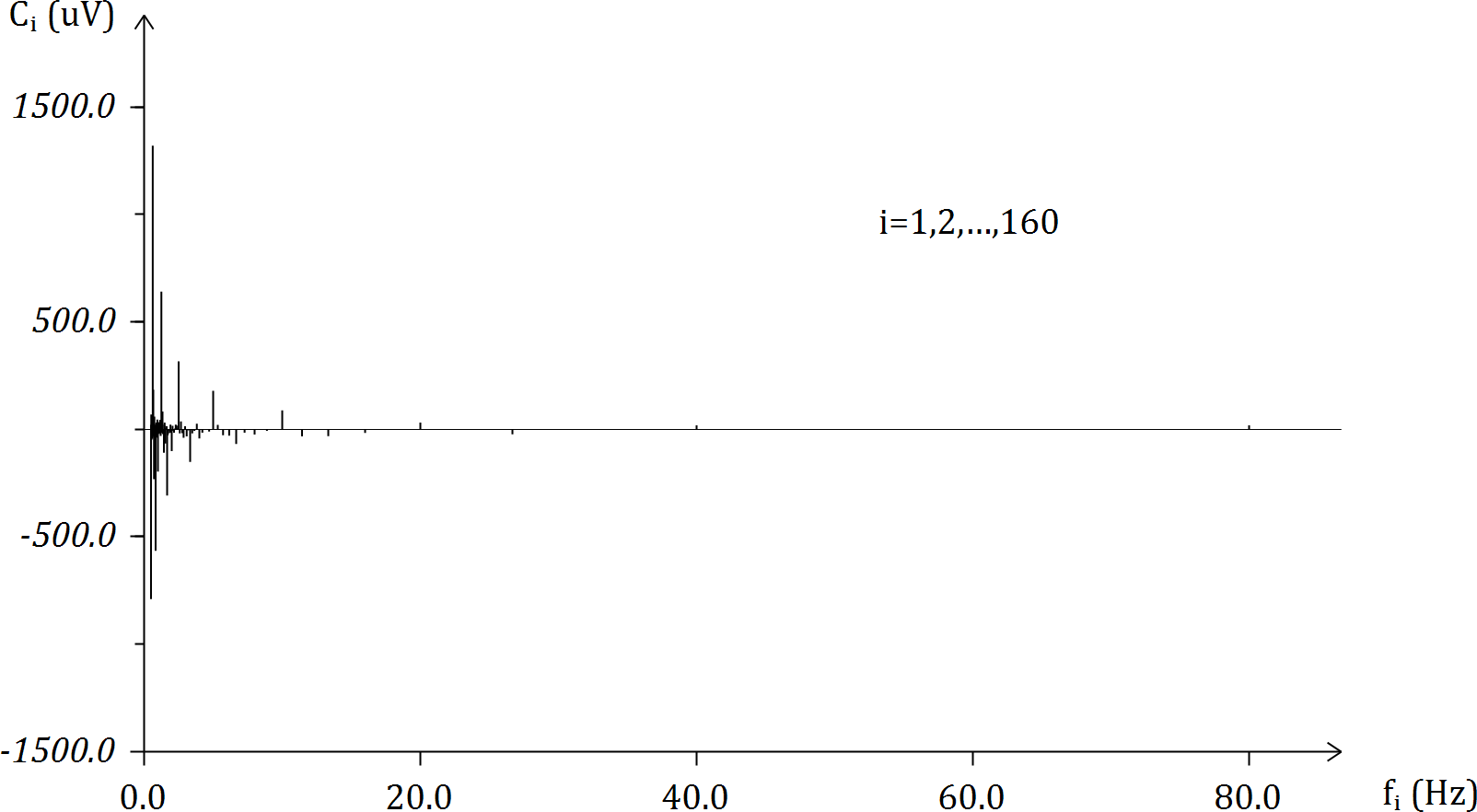}
\caption{The SWT corresponding to the sequence of ``samples'' represented in figure~\ref{f8}, from an electroencephalographic recording.}
\label{f9}
\end{figure}

In figure~\ref{f9}, the SWT is shown for the sequence of 160 ``samples'' (i.e., measured values) of that electroencephalographic recording.

\begin{figure}[H]
\centering
\includegraphics[width=4.5in]{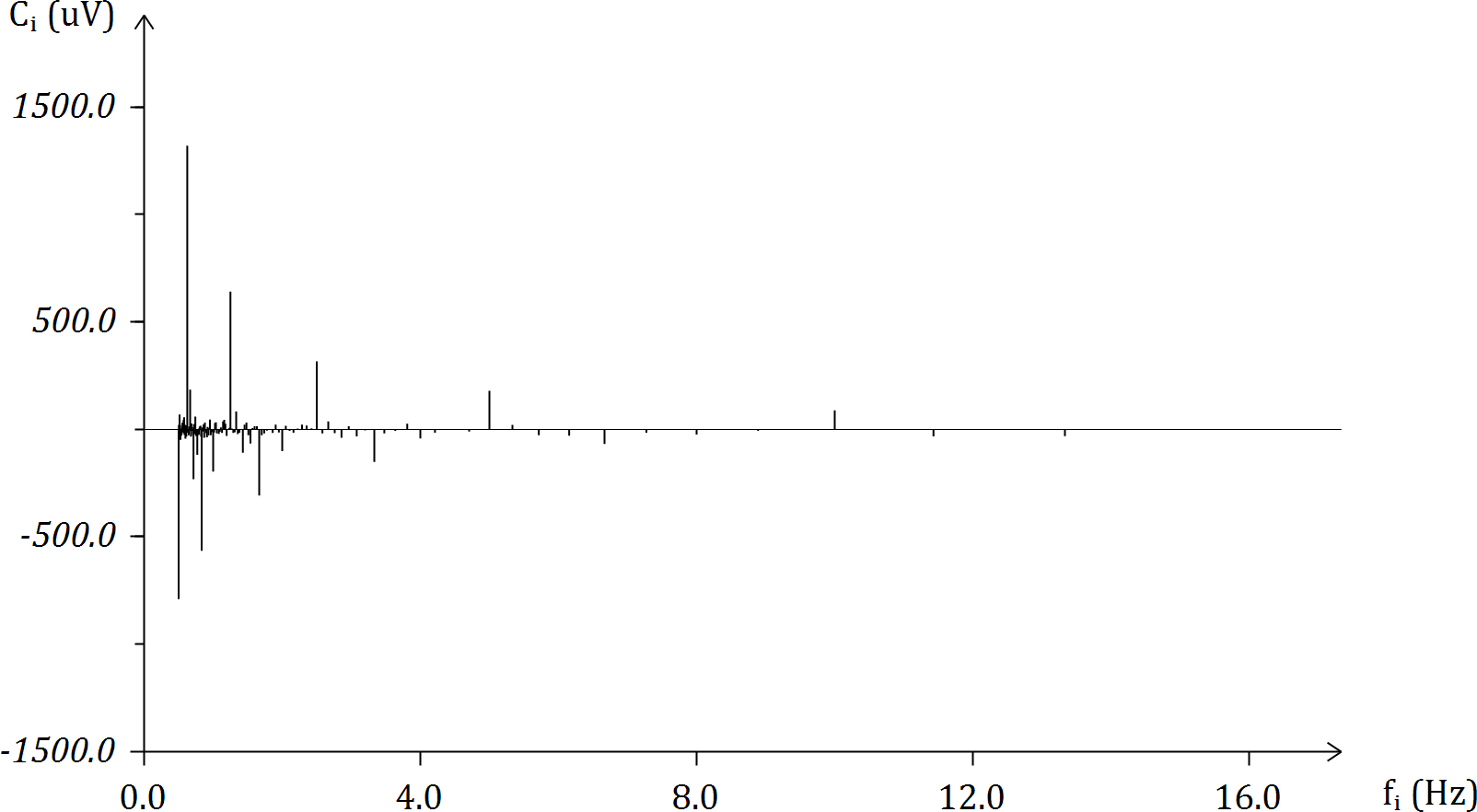}
\caption{Close-up of the SWT corresponding to the sequence of ``samples'' represented in figure~\ref{f9}, from an electroencephalographic recording, where $f_i<16$.}
\label{f11}
\end{figure}

In figure~\ref{f11}, a close-up of the most notable portion of the SWT is shown.

The sequence of 160 dyads such that the first element of the $i^{th}$ dyad ($i=1, 2, 3, \ldots, 160$) is $f_{i}$ (i.e., the frequency corresponding to $S_{i}$), and the second element of that dyad is $C_{i}$, is given below. (The 160 dyads are in increasing order according to the corresponding frequencies. These frequencies have been specified with a 4-digit accuracy.)

\vspace{4mm}

\begin{tabular}{ l l l l l }
$(0.5000; -789.5)$ & $(0.6667; 186.0)$ & $(1.0000; -195.5)$ & $(2.0000; -100.0)$ \\
$(0.5031; 21.5)$ & $(0.6723; 21.5)$ & $(1.0127; -13.5)$ & $(2.0513; 18.0)$ \\
$(0.5063; -11.0)$ & $(0.6780; -33.0)$ & $(1.0256; 30.5)$ & $(2.1053; -5.5)$ \\
$(0.5096; -22.0)$ & $(0.6838; 28.0)$ & $(1.0390; 33.5)$ & $(2.1622; -14.5)$ \\
$(0.5128; 70.0)$ & $(0.6897; -5.5)$ & $(1.0526; -18.0)$ & $(2.2222; 4.5)$ \\
$(0.5161; 20.5)$ & $(0.6957; 11.0)$ & $(1.0667; -1.0)$ & $(2.2857; 23.0)$ \\
$(0.5195; 20.5)$ & $(0.7018; -3.0)$ & $(1.0811; -19.5)$ & $(2.3529; 18.5)$ \\
$(0.5229; -2.0)$ & $(0.7080; -25.5)$ & $(1.0959; -10.0)$ & $(2.4242; 6.0)$ \\
$(0.5263; -47.5)$ & $(0.7143; -231.0)$ & $(1.1111; 10.0)$ & $(2.5000; 317.5)$ \\
$(0.5298; -26.5)$ & $(0.7207; 25.5)$ & $(1.1268; -16.0)$ & $(2.5806; -18.0)$ \\
$(0.5333; -19.0)$ & $(0.7273; 25.0)$ & $(1.1429; 38.5)$ & $(2.6667; 37.5)$ \\
$(0.5369; -28.5)$ & $(0.7339; -21.0)$ & $(1.1594; 46.0)$ & $(2.7586; -16.5)$ \\
$(0.5405; -23.0)$ & $(0.7407; 61.0)$ & $(1.1765; 27.0)$ & $(2.8571; -38.5)$ \\
$(0.5442; 5.5)$ & $(0.7477; -24.0)$ & $(1.1940; -29.0)$ & $(2.9630; 16.0)$ \\
$(0.5479; -15.5)$ & $(0.7547; -31.0)$ & $(1.2121; -0.5)$ & $(3.0769; -31.5)$ \\
$(0.5517; 12.0)$ & $(0.7619; -17.5)$ & $(1.2308; 8.5)$ & $(3.2000; -4.0)$ \\
$(0.5556; 33.0)$ & $(0.7692; -117.5)$ & $(1.2500; 642.0)$ & $(3.3333; -150.5)$ \\
$(0.5594; 29.0)$ & $(0.7767; -25.0)$ & $(1.2698; 5.0)$ & $(3.4783; -17.5)$ \\
$(0.5634; -12.5)$ & $(0.7843; -13.0)$ & $(1.2903; -15.5)$ & $(3.6364; -5.5)$ \\
$(0.5674; 28.5)$ & $(0.7921; -24.0)$ & $(1.3115; -12.5)$ & $(3.8095; 27.0)$ \\
$(0.5714; 45.5)$ & $(0.8000; 8.0)$ & $(1.3333; 84.5)$ & $(4.0000; -41.5)$ \\
$(0.5755; -17.0)$ & $(0.8081; 14.0)$ & $(1.3559; -21.0)$ & $(4.2105; -15.0)$ \\
$(0.5797; 57.5)$ & $(0.8163; 18.0)$ & $(1.3793; -14.5)$ & $(4.4444; -1.5)$ \\
$(0.5839; -25.0)$ & $(0.8247; -32.5)$ & $(1.4035; 3.0)$ & $(4.7059; -9.0)$ \\
$(0.5882; 22.5)$ & $(0.8333; -564.0)$ & $(1.4286; -107.5)$ & $(5.0000; 180.5)$ \\
$(0.5926; 2.5)$ & $(0.8421; 10.5)$ & $(1.4545; 23.5)$ & $(5.3333; 21.5)$ \\
$(0.5970; -42.0)$ & $(0.8511; -9.0)$ & $(1.4815; 32.5)$ & $(5.7143; -26.0)$ \\
$(0.6015; 19.0)$ & $(0.8602; 24.0)$ & $(1.5094; -26.0)$ & $(6.1538; -28.0)$ \\
$(0.6061; 20.0)$ & $(0.8696; -38.0)$ & $(1.5385; -64.5)$ & $(6.6667; -67.0)$ \\
$(0.6107; -32.5)$ & $(0.8791; 33.0)$ & $(1.5686; 6.5)$ & $(7.2727; -14.5)$ \\
$(0.6154; 5.5)$ & $(0.8889; -7.5)$ & $(1.6000; 15.5)$ & $(8.0000; -23.5)$ \\
$(0.6202; 18.0)$ & $(0.8989; -17.0)$ & $(1.6327; 15.5)$ & $(8.8889; -6.0)$ \\
$(0.6250; 1321.5)$ & $(0.9091; -36.0)$ & $(1.6667; -307.0)$ & $(10.0000; 89.0)$ \\
$(0.6299; -15.0)$ & $(0.9195; 12.5)$ & $(1.7021; -26.5)$ & $(11.4857; -31.5)$ \\
$(0.6349; -2.0)$ & $(0.9302; -27.0)$ & $(1.7391; -18.0)$ & $(13.3333; -31.0)$ \\
$(0.6400; 2.5)$ & $(0.9412; 7.0)$ & $(1.7778; -4.5)$ & $(16.0000; -15.5)$ \\
$(0.6452; -26.0)$ & $(0.9524; 46.5)$ & $(1.8182; 0.0)$ & $(20.0000; 32.5)$ \\
$(0.6504; 14.5)$ & $(0.9639; -26.5)$ & $(1.8605; -16.0)$ & $(26.6667; -22.5)$ \\
$(0.6557; -7.5)$ & $(0.9756; -12.0)$ & $(1.9048; 23.5)$ & $(40.0000; 19.5)$ \\
$(0.6612; 4.5)$ & $(0.9877; -9.0)$ & $(1.9512; -13.5)$ & $(80.0000; 20.0)$
\end{tabular}

\vspace{4mm}

\section{Discussion and prospects}

It must be emphasized that the SWTs of the corresponding approximations to a particular function have a pattern in common for high enough values of $N_s$. Although this topic will be addressed elsewhere, preliminary support for this will be given below.

Partial graphic representations of the SWTs displayed in figure~\ref{f7} above are shown in figure~\ref{f10}. The SWTs in which $N_s$ is equal to 100, 200 and 400, respectively, are partially presented in 11a, 11b, and 11c. In this case, the SWTs are described as ``partial'' because the axes of the abscissas extend only as far as the frequencies which are equal to or less than 15. To detect this pattern easily, the same scale has been used in the axes of the abscissas in 11a, 11b, and 11c.

\begin{figure}[H]
\centering
\subfloat[$N_s=100$.]{
\includegraphics[width=4.5in]{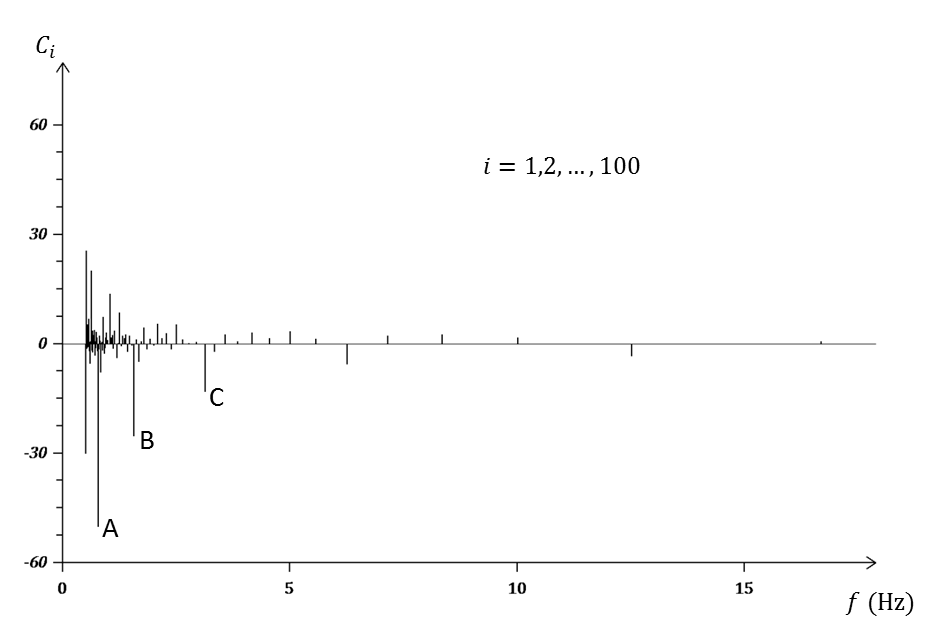}
\label{f10a}
}\qquad
\caption[]{}
\end{figure}
\begin{figure}[H]
\ContinuedFloat
\subfloat[$N_s=200$.]{
\includegraphics[width=4.5in]{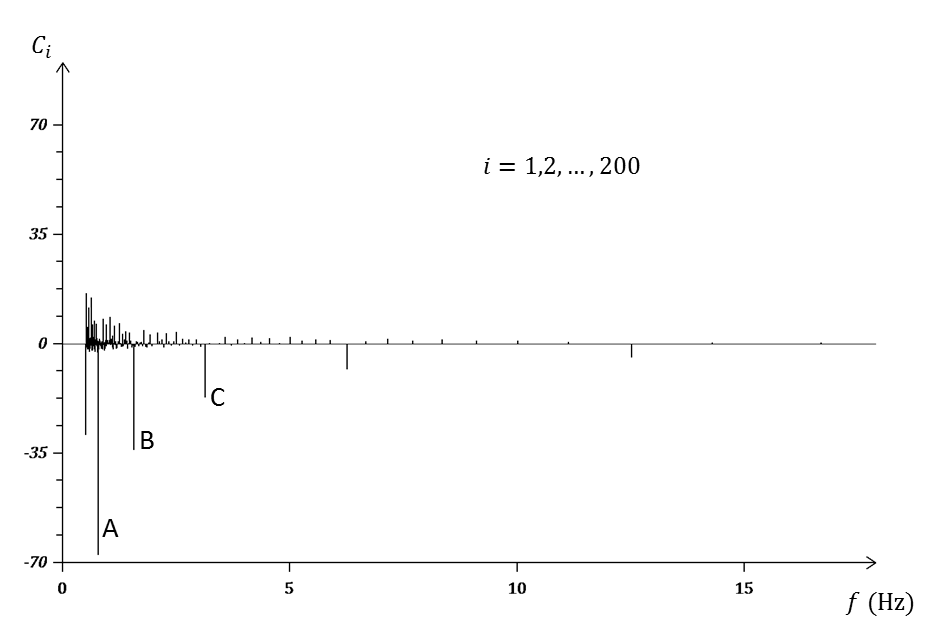}
\label{f10b}
}\qquad
\subfloat[$N_s=400$.]{
\includegraphics[width=4.5in]{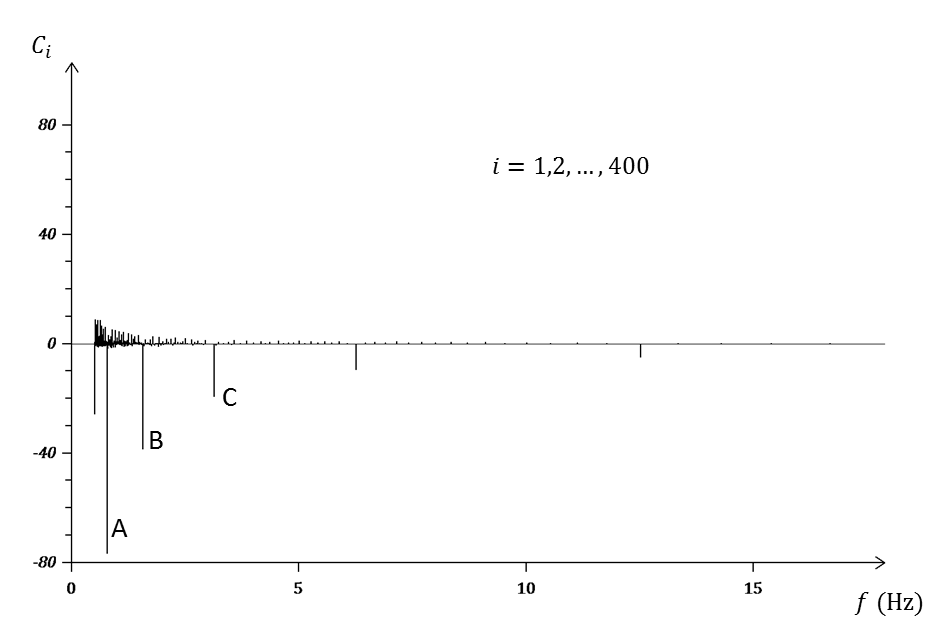}
\label{f10c}
}\qquad
\caption{The SWTs seen in figure~\ref{f7} are displayed partially in 11a, 11b, and 11c. The same scale has been used in the axes of the abscissas.}
\label{f10}
\end{figure}

Note, for example, the correspondence between the coefficients which have been indicated by the letter ``A'' in 11a, 11b and 11c. The correspondence between the coefficients indicated by ``B'' can also be seen, as can those indicated by the letter ``C". (Other interesting correspondences can also be observed, if desired.)

When comparing the SWTs corresponding to a given type of electroencephalographic recordings, care must be taken to use recordings made during the same $\Delta t$ and with the same sampling frequency.

The first of several computational tools for the use of the SWT to be made available for interested users has been installed on the website of the Applied Mathematics and Computer Simulation Group (www.appliedmath group.org). This tool is that of the Square Wave Transform (SWT) and makes it possible to obtain the SWT of electroencephalographic recordings automatically \cite{b6}.

%\subsection{A subsection}

%\newpage

\end{document}